\documentclass[11pt, fleqn]{amsart}

\usepackage[square,compress,comma, numbers,sort]{natbib}
\usepackage[colorlinks=true, citecolor=red, linkcolor=blue]{hyperref}
\usepackage{amsfonts,mathtools}
\usepackage{amsmath}
\usepackage{mathrsfs}

\allowdisplaybreaks[4]

\usepackage{amssymb}
\usepackage{color}
\usepackage{float}
\usepackage{eufrak}

\DeclarePairedDelimiterXPP\pk[1]{\mathbb{P}}\{ \}{}{ #1}
\DeclarePairedDelimiterXPP\E[1]{\mathbb{E}}\{ \}{}{	#1}

\usepackage{xparse}% http://ctan.org/pkg/xparse

\NewDocumentCommand{\ceil}{s O{} m}{%
	\IfBooleanTF{#1} % starred
	{\left\lceil#3\right\rceil} % \ceil*[..]{..}
	{#2\lceil#3#2\rceil} % \ceil[..]{..}
}
\NewDocumentCommand{\floor}{s O{} m}{%
	\IfBooleanTF{#1} % starred
	{\left\lfloor#3\right\rfloor}
	{#2\lfloor#3#2\rfloor}
}

\definecolor{c20}{rgb}{0.,0.7,0.}
\definecolor{c30}{rgb}{0.,0.,1.}
\definecolor{c40}{rgb}{1,0.1,0.7}
\definecolor{c50}{rgb}{1,0,0}
\definecolor{c60}{rgb}{1,0.9,0.1}
\definecolor{c70}{rgb}{0.50,1.00,0.00}

\numberwithin{equation}{section}
\newtheorem{theo}{Theorem}[section]
\newtheorem{sat}[theo]{Proposition}
\newtheorem{de}[theo]{Definition}
\newtheorem{lem}[theo]{Lemma}

\newtheorem{example}[theo]{Example}
\newtheorem{korr}[theo]{Corollary}
\newtheorem{remark}[theo]{Remark}

\numberwithin{equation}{section}

\newcommand{\prooftheo}[1]{ \textsc{Proof of Theorem} \ref{#1} }

\newcommand{\prooflem}[1]{\textsc{Proof of Lemma} \ref{#1}}

\newcommand{\QED}{\hfill $\Box$}

\newcommand{\COM}[1]{}

\def\IF{\infty}

\newcommand{\R}{\mathbb{R}}

\newcommand{\BQN}{\begin{eqnarray}}
\newcommand{\EQN}{\end{eqnarray}}
\newcommand{\BQNY}{\begin{eqnarray*}}
	\newcommand{\EQNY}{\end{eqnarray*}}

\newcommand{\limit}[1]{\lim_{#1 \to   \infty}}

\newcommand{\kb}[1]{\boldsymbol{#1}}
\newcommand{\vk}[1]{\kb{#1}}

\def\bqny#1{{\begin{eqnarray*} #1 \end{eqnarray*}}}
\def\bqn#1{{\begin{eqnarray} #1 \end{eqnarray}}}

\newcommand{\BS}{\begin{sat}}
	\newcommand{\ES}{\end{sat}}
\newcommand{\BT}{\begin{theo}}
	\newcommand{\ET}{\end{theo}}
\newcommand{\BK}{\begin{korr}}
	\newcommand{\EK}{\end{korr}}

\newcommand{\BEX}{\begin{example}}
	\newcommand{\EEX}{\end{example}}

\newcommand{\BD}{\begin{de}}
	\newcommand{\ED}{\end{de}}

\newcommand{\BIT}{\begin{itemize}}
	\newcommand{\EIT}{\end{itemize}}
\newcommand{\BDI}{\begin{description}}
	\newcommand{\EDI}{\end{description}}

\newcommand{\BRM}{\begin{remark}}
	\newcommand{\ERM}{\end{remark}}

\newcommand{\BEL}{\begin{lem}}
	\newcommand{\EEL}{\end{lem}}

\newcommand{\nelem}[1]{{Lemma \ref{#1}}}

\begin{document}

\title{Two-dimensional Brownian risk model for cumulative Parisian ruin probability}

\author{Konrad Krystecki}
\address{Konrad Krystecki, Department of Actuarial Science, %\\Faculty of Business and Economics\\
	University of Lausanne,\\
	UNIL-Dorigny, 1015 Lausanne, Switzerland and
	 Mathematical Institute, University of Wroc\l aw, pl. Grunwaldzki 2/4, 50-384 Wroc\l aw, Poland
}
\email{Konrad.Krystecki@unil.ch}

\bigskip

\date{\today}
 \maketitle

 {\bf Abstract:} Let $(W_1(s), W_2(t)), s,t\ge 0$ be a bivariate Brownian motion with standard Brownian motion marginals and constant correlation $\rho \in (-1,1).$ In this contribution we derive precise approximations for cumulative Parisian ruin conditioned on the occurrence of the ruin of the aforementioned two-dimensional Brownian motion, i.e.
 $$\pk*{\begin{array}{ccc}\int_{[0,1]} \mathbf{1}(W_1^*(s)>u)ds>H_1(u) \\ \int_{[0,1]} \mathbf{1}(W_2^*(t)>au)dt>H_2(u)\end{array}\Bigg{|}\exists_{v,w \in [0,1]}\begin{array}{ccc} W_1(v)-c_1v>u \\ W_2(w)-c_2w>au \end{array}}.$$
 We study the asymptotics for specific functions $\vk H(u)$ for $u$ being proportional to initial position of the Brownian motion, which determines how long does the process need to spend over the barrier.



\section{Introduction}
Consider the following Brownian risk model for two portfolios
$$R_i(t)=u_i+c_i t - W_i(t), i=1,2, $$
where the claims $W_i(t), t \ge 0$ are modeled by two dependent standard Brownian motions, initial capitals $u_i > 0$ and premium rates $c_i.$ We model dependence between coordinates as it was introduced in e.g. \cite{delsing2018asymptotics} and \cite{DIEKER2005}
\BQNY \label{BB}
(W_1(s),W_2(t))=(B_1(s), \rho B_1(t)+ \sqrt{1- \rho^2} B_2(t)), \quad s,t\ge 0,
\EQNY
where $B_1,B_2$ are two independent standard Brownian motions and $\rho \in [-1,1]$. The probability of ruin of a single portfolio in the finite time horizon is given by (see e.g., \cite{MandjesKrzys})
\BQNY \label{single}
{ \pi}_T(c_i,u):=\pk*{ \inf_{t\in [0,T]} R_i(t) < 0}&=& \pk*{\sup_{t\in [0,T]} (W_i(t)- c_i t)> u}\notag \\
&=& \Phi\left(-\frac{ u}{ \sqrt{T}} -{c_i\sqrt{T}}\right)+
e^{-2c_iu}\Phi\left(- \frac{ u}{\sqrt{T}} +{c_i\sqrt{T}}\right)
\EQNY
for  $i=1,2, u\geq 0$, with $\Phi$ the distribution function of an $N(0,1)$ random variable. Since from self-similarity of Brownian motion we have that for $c_1'=\frac{c_1}{\sqrt{T}}, u'=\frac{u}{\sqrt{T}}$
$$B(tT)-c_1 t >u \stackrel{D}{\Leftrightarrow} \sqrt{T}B(t) - c_1 t > u \Leftrightarrow  B(t) - c_1' t > u',$$
then without loss of generality one can assume $T=1.$ Denote $W_i^*(s)= W_i(s)-c_i s, i=1,2.$ In the literature several models describing the ruin have been introduced and investigated for the multidimensional setting. For example define the simultaneous ruin probability as
$$ \overline\pi_{A,\rho}(c_1,c_2,u, v)= \pk*{ \exists_{ s \in A}:  W_1^*(s) > u,  W_2^*(s) > v}$$
which has been recently studied in \cite{SIM} for $A=[0,1]$. Similarly, define non-simultaneous ruin probability as
$$ \pi_{A \times B, \rho}(c_1,c_2,u,v)= \pk*{ \exists_{ s \in A, t \in B}:  W_1^*(s) > u,  W_2^*(t) > v}$$
which has been studied for the case $A=B=[0,1]$ in \cite{DHK20}.

In this contribution we focus on investigating different definition of the non-simultaneous ruin for two-dimensional risk portfolios previously investigated in e.g. \cite{KEZX17}. Let us introduce a ruin definition that generalises the notion of the Parisian ruin (see e.g. \cite{loeffen2013,krystecki2021parisian}) and since it allows the time above the threshold to come from disjoint intervals, it is also called cumulative Parisian ruin. For one dimensional model it can be formulated in the following way
\BQNY \label{singleSojurn}
\mathcal{S}_{A,H(u)}(c,u) := \pk*{\int_{A} \mathbf{1}(W^*(s)>u)ds>H(u)},
\EQNY
for some $H(u) \ge 0$ and $A=[0,T].$ Similarly, the two-dimensional model follows
\BQNY \label{doubleSojurn}
\mathcal{S}_{A \times B,\vk H(u)}(c_1,c_2,u,au) := \pk*{\int_{A} \mathbf{1}(W_1^*(s)>u)ds>H_1(u),\int_{B} \mathbf{1}(W_2^*(t)>au)dt>H_2(u)},
\EQNY
for some $H_1(u),H_2(u) \ge 0$ and compact sets $A, B.$ Clearly, we have that cumulative Parisian ruin is a generalisation of classical ruin since
$$\pi_{[0,1]^2,\rho}(c_1,c_2,u,au)=\mathcal{S}_{[0,1]^2,0}(c_1,c_2,u,au). $$
For more general compact sets $A,B$ we have the following relation
\BQN \label{comparison}
\pi_{A \times B,\rho}(c_1,c_2,u,au) \ge \mathcal{S}_{A \times B,\vk H(u)}(c_1,c_2,u,au).
\EQN
In this contribution we aim at understanding the exact relation between different types of ruin through the asymptotic behaviour of
\BQNY
\lefteqn{\mathscr{S}^*_{[0,1]^2,\vk H(u)}(c_1,c_2,u,au):=}&&\\
&&\pk*{\int_{[0,1]} \mathbf{1}(W_1^*(s)>u)ds>H_1(u),\int_{[0,1]} \mathbf{1}(W_2^*(t)>au)dt>H_2(u)\Bigg{|}\exists_{v,w \in [0,1]}\begin{array}{ccc} W_1^*(v)>u \\ W_2^*(w)>au \end{array}}.
\EQNY
In the one-dimensional context it was proven that the cumulative Parisian ruin differs from classical ruin only in constant. For two-dimensional case we can similarly write that
\BQNY
\mathscr{S}_{S_1,S_2}(c_1,c_2,u,au) &=& \frac{\mathcal{S}_{S_1,S_2}(c_1,c_2,u,au)}{\pi_{[0,1]^2, \rho}(c_1,c_2,u,au)}
\EQNY
and by proving that $\mathscr{S}_{S_1,S_2}(c_1,c_2,u,au)$ is constant we prove that the cumulative Parisian ruin and classical ruin are of the same order. From \cite{dkebicki2016parisian} and \cite{ji2018cumulative} we have that usually Parisian ruin differs from cumulative Parisian ruin only in constant and in \cite{krystecki2021parisian} we have seen that for non-simultaneous model Parisian ruin differs from classical ruin only in constant. Arbitrary choice of $H(u) = \frac{(S_1,S_2)}{u^2}$ is closely connected to the behaviour of variance for the Brownian motion (see \cite{pickands1969}). For the choice of $H(u)=o(\frac{1}{u^2})$ following the same line of proof we have that
$$\lim_{u\to \IF}\mathscr{S}^*_{[0,1]^2,\vk H(u)}(c_1,c_2,u,au)=1,$$
which is a natural result, since the required period of crossing the barrier is so short that it is enough that the process crosses it once to stay there for the required time. On the other hand, if we choose $H(u)$ such that $\frac{H(u)}{\frac{1}{u^2}} \rightarrow \IF, \forall_u H(u)<1,$ then the methods employed in this contribution are not sufficient and the asymptotics are of different order.

\section{Main results}
We begin with dimension-reduction cases, where one of the coordinates dominates the other and the results are up to constant the same as one-dimensional results. This behavior was already observed in \cite{DHK20} for the classical ruin and in \cite{krystecki2021parisian} for Parisian ruin.

\BT \label{simpleSojurn} If $a \le \rho, $ then
\bqny{
\lim_{u \to \IF}\mathscr{S}_{S_1,S_2}(c_1,c_2,u,au) &=& \frac{((2+S_1)\Phi(\sqrt{\frac{S_1}{2}})-\sqrt{\frac{S_1}{\pi}}e^{-\frac{S_1}{4}})}{2}.
}
\ET
We define constants for cumulative Parisian ruin as follows
$$\widehat{\mathcal{P}}(w_1,w_2, f(u) ):=\int_{\R}\pk*{ \int_{[0,\infty)}\mathbf{1}( B(s) - w_1s> x)ds>f(u)} e^{ w_2x} dx, $$
$$\widehat{\mathcal{H}}(w_1,w_2, f(u) ):=\lim_{\Delta \to \infty}\int_{\R}\frac{1}{\Delta}\pk*{ \int_{[0,\Delta]}\mathbf{1}( B(t) - w_1t> x)dt>f(u)} e^{ w_2x} dx,$$
$$\widehat{\mathcal{R}}(S_1,S_2)=\int_{\R^2}\pk*{\begin{array}{ccc}\int_{[0,\infty)}\mathbf{1}(W_1(s)-s>x)ds>S_1 \\ \int_{[0,\infty)}\mathbf{1}(W_2(t)-at> y)dt>S_2\end{array}} e^{\frac{1-a\rho}{1-\rho^2} x + \frac{a-\rho}{1-\rho^2} y} dxdy \in (0,\IF). $$
In each particular case, finintess and positivity of $\widehat{\mathcal{P}}$ and $\widehat{\mathcal{H}}$ has been proven in \nelem{1dFiniteSojurn}.
\BT \label{MainSojurn} Let   $\rho \in (-1,1)$  and $a\in  (\max(0,\rho), 1]$ be given. \\
(i) If $\rho> A_a $, then
\bqn{ \label{M1}
	\lim_{u \to \IF}\mathscr{S}_{S_1,S_2}(c_1,c_2;u,au) = \frac{\widehat{\mathcal{R}}_{S_1,S_2}}{\widehat{\mathcal{R}}_{0,0}}.
}
\newline
(ii) If $\rho= A_a $ and $a <1$, then
\bqn{ \label{M2} \lim_{u \to \IF}\mathscr{S}_{S_1,S_2}(c_1,c_2;u,au) = \frac{(1-a\rho)\widehat{\mathcal{P}}(\frac{1-a\rho}{1-\rho^2},\frac{1-a\rho}{1-\rho^2},S_1)\widehat{\mathcal{H}}(a,2a,S_2)}{2a(1-\rho^2)}.
}
(iii) If $\rho=A_a, a=1$, then
\bqn{ \label{M4}
	\lim_{u \to \IF}\mathscr{S}_{S_1,S_2}(c_1,c_2;u,au) =  \frac{C_{3,1} C'_{3,1}+C_{3,2} C'_{3,2}}{C_3},
}
%%%%%%%trzeba rozdzielić stałe, bo S_1,S_2
where $C_{3,1}=\widehat{\mathcal{P}}(2,2,S_1) \widehat{\mathcal{H}}(1,2,S_2),C_{3,2}=\widehat{\mathcal{P}}(2,2,S_2) \widehat{\mathcal{H}}(1,2,S_1)$ and
$$C'_{3,1}=
\begin{cases}
e^{-2\frac{(\frac{1}{2}c_1+c_2)^2}{3}} \Phi\left(c_2+\frac{1}{2}c_1\right), &-\frac{1}{2}c_1<c_2 \\
1, & otherwise,
\end{cases}, \quad C'_{3,2}=
\begin{cases}
e^{-2\frac{(\frac{1}{2}c_2+c_1)^2}{3}}  \Phi\left(c_1+\frac{1}{2}c_2\right), & -\frac{1}{2}c_2<c_1\\
1, & otherwise,
\end{cases} $$
$$ C_3=
\begin{cases}
e^{-2\frac{(\frac{1}{2}c_1+c_2)^2}{3}} \Phi\left(c_2+\frac{1}{2}c_1\right)
+e^{-2\frac{(\frac{1}{2}c_2+c_1)^2}{3}}  \Phi\left(c_1+\frac{1}{2}c_2\right), & c_2>\max(-\frac{1}{2}c_1,-2c_1)\\
e^{-2\frac{(\frac{1}{2}c_1+c_2)^2}{3}} \Phi\left(c_2+\frac{1}{2}c_1\right)
+\frac{1}{2}, &-\frac{1}{2}c_1<c_2 \le -2c_1 \\
\frac{1}{2}
+e^{-2\frac{(\frac{1}{2}c_2+c_1)^2}{3}}  \Phi\left(c_1+\frac{1}{2}c_2\right), &-2c_1<c_2 \le -\frac{1}{2}c_1 \\
1 , &c_2\le \min(-\frac{1}{2}c_1,-2c_1).
\end{cases} $$
(iv) If $a<1, \rho< A_a $, then
\bqn{\label{M5} \lim_{u \to \IF}\mathscr{S}_{S_1,S_2}(c_1,c_2;u,au) = -\frac{\widehat{\mathcal{P}}(\frac{1-a\rho}{1-\rho^2 t_*},\frac{1-a\rho}{1-\rho^2 t_*},S_1) \widehat{\mathcal{H}}(\frac{a}{t_*}, \frac{2a}{t_*},S_2)}{2\rho}.}
(v) If $a=1, \rho< A_a $, then
\bqn{\label{M6} \lim_{u \to \IF}\mathscr{S}_{S_1,S_2}(c_1,c_2;u,au) =  -\frac{C_5}{2\rho},} where
 $t_*=\frac{1}{\rho(2\rho-1)}, C_5=\begin{cases} \widehat{\mathcal{P}}(\frac{1-\rho}{1-\rho^2 t_*},\frac{1-\rho}{1-\rho^2 t_*},S_1) \widehat{\mathcal{H}}(\frac{1}{t_*}, \frac{2}{t_*},S_2)&c_1 \le c_2 \\ \widehat{\mathcal{P}}(\frac{1-\rho}{1-\rho^2 t_*},\frac{1-\rho}{1-\rho^2 t_*},S_2) \widehat{\mathcal{H}}(\frac{1}{t_*}, \frac{2}{t_*},S_1),& c_1 > c_2 \end{cases}.$
\ET
\section{Proofs}
We begin with the proofs of dimension-reduction cases, where we prove that one of the coordinates only contributes to a constant. Notice that
\BQNY
\mathscr{S}_{S_1,S_2}(c_1,c_2,u,au) &=& \frac{\mathcal{S}_{S_1,S_2}(c_1,c_2,u,au)}{\pi_{[0,1]^2, \rho}(c_1,c_2,u,au)}.
\EQNY
Since the asymptotics of $\pi_{[0,1]^2, \rho}(c_1,c_2,u,au)$ has already been studied in \cite{DHK20} it is often easier to analyze the asymptotics of $\mathcal{S}_{S_1,S_2}(c_1,c_2,u,au),$ hence in the proofs we focus on investigating the behaviour of the latter.
\subsection{Proof of Theorem \ref{simpleSojurn}}
\underline{Case (i): $a < \rho.$} Notice that
$$\mathcal{S}_{[0,1]^2,\frac{(S_1,S_2)}{u^2}}(c_1,c_2,u,au) \le  \mathcal{S}_{[0,1],\frac{S_1}{u^2}}(c_1,u).$$
Further we have for large enough $u$
\BQNY
\lefteqn{\mathcal{S}_{[0,1]^2,\frac{S}{u^2}}(c_1,c_2,u,au)} && \\ &\ge& \pk*{\begin{array}{ccc}\int_{[0,1]} \mathbf{1}(B_1^*(s)>u,\rho B_1(s) + \sqrt{1-\rho^2} B_2(s) - c_2 s>au)ds>\frac{S_1}{u^2}\\\int_{[0,1]} \mathbf{1}(B_1^*(s)>u-\frac{1}{\sqrt{u}},\rho B_1(s) + \sqrt{1-\rho^2} B_2(s) - c_2 s>au)ds>\frac{\max(S_1,S_2)}{u^2} \end{array}}\\
&\ge&\pk*{\begin{array}{ccc}\int_{[0,1]} \mathbf{1}(B_1^*(s)>u,\rho (u+c_1s) + \sqrt{1-\rho^2} B_2(s) - c_2 s>au)ds>\frac{S_1}{u^2} \\ \int_{[0,1]} \mathbf{1}(B_1^*(s)>u-\frac{1}{\sqrt{u}},\rho (u-\frac{1}{\sqrt{u}}+c_1s) + \sqrt{1-\rho^2} B_2(s) - c_2 s>au)ds>\frac{\max(S_1,S_2)}{u^2}\end{array}}\\
&\ge&\pk*{\begin{array}{ccc}\int_{[0,1]} \mathbf{1}(B_1^*(s)>u)\mathbf{1}\left(\forall_{t \in [0,1]} B_2(t)>\frac{(a-\rho)u+(c_2-\rho c_1)t}{\sqrt{1-\rho^2}}\right)ds>\frac{S_1}{u^2}\\\int_{[0,1]} \mathbf{1}(B_1^*(s)>u-\frac{1}{\sqrt{u}})\mathbf{1}\left(\forall_{t \in [0,1]} B_2(t)>\frac{(a-\rho)u+(c_2-\rho c_1)t}{\sqrt{1-\rho^2}}\right)ds>\frac{\max(S_1,S_2)}{u^2} \end{array}}\\
&=&\pk*{\begin{array}{ccc}\mathbf{1}\left(\forall_{t \in [0,1]}B_2(t)>\frac{(a-\rho)u+(c_2-\rho c_1)t}{\sqrt{1-\rho^2} }\right)\int_{[0,1]} \mathbf{1}(B_1^*(s)>u)ds>\frac{S_1}{u^2}\\\mathbf{1}\left(\forall_{t \in [0,1]}B_2(t)>\frac{(a-\rho)u+(c_2-\rho c_1)t}{\sqrt{1-\rho^2} }\right)\int_{[0,1]} \mathbf{1}(B_1^*(s)>u-\frac{1}{\sqrt{u}})ds>\frac{\max(S_1,S_2)}{u^2} \end{array}}\\
&=&\pk*{\forall_{t \in [0,1]}B_2(t)>\frac{(a-\rho)u+(c_2-\rho c_1)t}{\sqrt{1-\rho^2} }}\pk*{\begin{array}{ccc}\int_{[0,1]} \mathbf{1}(B_1^*(s)>u)ds>\frac{S_1}{u^2} \\ \int_{[0,1]} \mathbf{1}(B_1^*(s)>u-\frac{1}{\sqrt{u}})ds>\frac{\max(S_1,S_2)}{u^2} \end{array}}.
\EQNY
Since $a<\rho$, we have that
$$\lim_{u \to \IF}\pk*{\forall_{t \in [0,1]}B_2(t)>\frac{(a-\rho)u+(c_2-\rho c_1)t}{\sqrt{1-\rho^2} }}=1.$$
Further by the independence of increments and self-similarity of Brownian motion we have that for $\widehat{B_1}$ a standard Brownian motion independent of $B_1$ and for $t_{\inf}=\inf\{s \in [0,1]:B_1^*(s)=u|\exists_{t \in [0,1]}:B_1^*(t)>u\}$
\BQN \label{ineqs}
\lefteqn{\pk*{\begin{array}{ccc} \int_{[0,1]} \mathbf{1}(B_1^*(s)>u)ds>\frac{S_1}{u^2} \\ \int_{[0,1]} \mathbf{1}(B_1^*(s)>u-\frac{1}{\sqrt{u}})ds>\frac{\max(S_1,S_2)}{u^2}\end{array}}}&& \nonumber\\
&=& \pk*{\begin{array}{ccc} \int_{[t_{\inf},1]} \mathbf{1}(B_1^*(s)>u)ds>\frac{S_1}{u^2} \\ \int_{[0,1]} \mathbf{1}(B_1^*(s)>u-\frac{1}{\sqrt{u}})ds>\frac{\max(S_1,S_2)}{u^2}\end{array}} \nonumber \\
&\ge& \pk*{\begin{array}{ccc} \int_{[t_{\inf},1]} \mathbf{1}(B_1^*(s)-B_1^*(t_{\inf})+B_1^*(t_{\inf})>u)ds>\frac{S_1}{u^2} \\ \forall_{s \in (t_{\inf}-\frac{\max(0,S_2-S_1)}{u^2},t_{\inf})}B_1^*(s)-B_1^*(t_{\inf})+B_1^*(t_{\inf})>u-\frac{1}{\sqrt{u}} \end{array}} \nonumber\\
&=& \pk*{\begin{array}{ccc} \int_{[t_{\inf},1]} \mathbf{1}(B_1^*(s)-B_1^*(t_{\inf})>0)ds>\frac{S_1}{u^2} \\ \forall_{s \in (t_{\inf}-\frac{\max(0,S_2-S_1)}{u^2},t_{\inf})}B_1^*(t_{\inf})-B_1^*(s)<\frac{1}{\sqrt{u}} \end{array}} \nonumber\\
&=& \pk*{\int_{[t_{\inf},1]} \mathbf{1}(B_1^*(s)-B_1^*(t_{\inf})>0)ds>\frac{S_1}{u^2}}\pk*{\forall_{s \in (t_{\inf}-\frac{\max(0,S_2-S_1)}{u^2},t_{\inf})}B_1^*(t_{\inf})-B_1^*(s)<\frac{1}{\sqrt{u}}} \nonumber\\
&=& \pk*{\int_{[t_{\inf},1]} \mathbf{1}(B_1^*(s)-B_1^*(t_{\inf})+B_1^*(t_{\inf})>u)ds>\frac{S_1}{u^2}}\pk*{\forall_{s \in (0,\max(\frac{S_2-S_1}{u^2},0))} \widehat{B_1}(s)+c_1s<\frac{1}{\sqrt{u}}} \nonumber\\
&=& \pk*{\int_{[0,1]} \mathbf{1}(B_1^*(s)>u)ds>\frac{S_1}{u^2}}\pk*{\forall_{s \in [0,\max(S_2-S_1,0)]} \widehat{B_1}(s)+\frac{c_1s}{u}<\sqrt{u}}.
\EQN
Finally we have that
\BQNY
\lim_{u \to \IF}\pk*{\forall_{s \in [0,\max(S_2-S_1,0)]} B(s)+\frac{c_1s}{u}<\sqrt{u}}&=&1.
\EQNY
One dimensional asymptotics found in \cite{debicki2019sojourn}[Cor 3.2, Thm 3.7] gives us
\BQN \label{1dressoj}
\pk*{\int_{[0,1]} \mathbf{1}(B_1^*(s)>u)ds>\frac{S_1}{u^2}} &\sim& \int\limits_{\mathbb{R}}e^x\pk{\int\limits_{0}^{\IF}\mathbf{1}(B(t)-t>x)dt>S_1}dx \Psi(u+c_1)\nonumber \\
 &\sim& ((2+S_1)\Phi(\sqrt{\frac{S_1}{2}})-\sqrt{\frac{S_1}{\pi}}e^{-\frac{S_1}{4}})\Psi(u+c_1)
\EQN
and hence the proof of case (i) is complete.
\newline
\underline{Case (ii): $a = \rho.$} For $\Delta>0$ we have
\BQNY \mathcal{S}_{[0,1]^2,\frac{S_1,S_2}{u^2}}(c_1,c_2;u, au) &\le& \mathcal{S}_{[1-\frac{\Delta}{u^2},1]^2,\frac{S_1}{u^2}}(c_1,c_2;u, au)+ \pi_{[0,1]^2 \setminus [1-\frac{\Delta}{u^2},1]^2,\rho}(c_1,c_2;u, au)\\
& := & \mathbb{P}+\pi_{[0,1]^2 \setminus [1-\frac{\Delta}{u^2},1]^2,\rho}(c_1,c_2;u, au).
\EQNY
Denote $\overline{\Delta}(u)=[1-\frac{1}{\sqrt{u}},1].$ Then
\BQNY \mathbb{S} &\le& \pk*{\int\limits_{\overline{\Delta}(u)} \mathbf{1}(W_1^*(s)>u)ds>\frac{S_1}{u^2},\int\limits_{\overline{\Delta}(u)} \mathbf{1}(W_2^*(t)>au)dt>\frac{S_2}{u^2}, \forall_{v \in \overline{\Delta}(u)} W_1^*(v)<u+\frac{1}{\sqrt{u}}}\\
&+& \pk*{\exists_{v \in \overline{\Delta}(u)} W_1^*(v)>u+\frac{1}{\sqrt{u}}}:=\mathbb{P}_1+\mathbb{P}_2.
\EQNY
Further
\BQNY
\mathbb{P}_1&\le&\pk*{\begin{array}{ccc}\int\limits_{\overline{\Delta}(u)} \mathbf{1}(W_1^*(s)>u)ds>\frac{S_1}{u^2} \\ \int\limits_{\overline{\Delta}(u)} \mathbf{1}(B_2(t)-\frac{c_2-\rho c_1}{1-\rho^2}t>\frac{-\rho}{\sqrt{u}})dt>\frac{S_2}{u^2}\end{array},\forall_{v \in \overline{\Delta}(u)} W_1^*(v)<u+\frac{1}{\sqrt{u}}}\\
&\le&\pk*{\int\limits_{\overline{\Delta}(u)} \mathbf{1}(W_1^*(s)>u)ds>\frac{S_1}{u^2}}\pk*{\int\limits_{\overline{\Delta}(u)} \mathbf{1}(B_2(t)-\frac{c_2-\rho c_1}{1-\rho^2}t>\frac{-\rho}{\sqrt{u}})dt>\frac{S_2}{u^2}}\\
&\le&\pk*{\int\limits_{\overline{\Delta}(u)} \mathbf{1}(W_1^*(s)>u)ds>\frac{S_1}{u^2}}\pk*{\sup_{s \in \overline{\Delta}(u)} B_2(t)-\frac{c_2-\rho c_1}{1-\rho^2}t>\frac{-\rho}{\sqrt{u}}}\\
&\sim&\pk*{\int\limits_{\overline{\Delta}(u)} \mathbf{1}(W_1^*(s)>u)ds>\frac{S_1}{u^2}}\pk*{\sup_{s \in \overline{\Delta}(u)} B_2(t)-\frac{c_2-\rho c_1}{1-\rho^2}t>0}\\
&=&\mathcal{S}_{[0,1],\frac{S_1}{u^2}}(c_1,u)\Phi\left( \frac{\rho c_1-c_2 }{\sqrt{1-\rho^2}}\right)(1+o(1)), \quad u \to \IF.\\
\EQNY
From \cite{DHK20} [Thm 2.1] we have that
$$\mathbb{P}_2=o\left(\pi_{\overline{\Delta}(u),\rho}(c_1,u)\right).$$
Further with  \cite{debicki2019sojourn}[Thm 3.7] we have that for some $C>0$
$$\lim_{u \to \IF}\frac{\mathcal{S}_{\overline{\Delta}(u),\frac{S_1}{u^2}}(c_1,u)}{\pi_{\overline{\Delta}(u),\rho}(c_1,u)}=C$$
and hence as $u \to \IF$
\BQNY \mathcal{S}_{[0,1]^2,\frac{S_1,S_2}{u^2}}(c_1,c_2;u, au) &\le& \Phi\left( \frac{\rho c_1-c_2 }{\sqrt{1-\rho^2}}\right)\mathcal{S}_{[0,1],\frac{S_1}{u^2}}(c_1,u).
\EQNY
For the lower bound notice that similarly as in case $a<\rho$
\BQNY \mathcal{S}_{[0,1]^2,\frac{S}{u^2}}(c_1,c_2;u, au) &\ge& \pk*{\forall_{t \in \overline{\Delta}(u)}B_2(t)>\frac{(c_2-\rho c_1)t}{\sqrt{1-\rho^2} }}\\
&&\times \pk*{\begin{array}{ccc}\int_{\overline{\Delta}(u)} \mathbf{1}(B_1^*(s)>u)ds>\frac{S_1}{u^2} \\ \int_{\overline{\Delta}(u)} \mathbf{1}(B_1^*(s)>u-\frac{1}{\sqrt{u}})ds>\frac{\max(S_1,S_2)}{u^2} \end{array}}.
\EQNY
With \eqref{ineqs} we have
$$\pk*{\begin{array}{ccc}\int_{\overline{\Delta}(u)} \mathbf{1}(B_1^*(s)>u)ds>\frac{S_1}{u^2} \\ \int_{\overline{\Delta}(u)} \mathbf{1}(B_1^*(s)>u-\frac{1}{\sqrt{u}})ds>\frac{\max(S_1,S_2)}{u^2} \end{array}} \sim \mathcal{S}_{[0,1],\frac{S_1}{u^2}}(c_1,u)$$
and further we notice that
$$\pk*{\forall_{t \in \overline{\Delta}(u)}B_2(t)>\frac{(c_2-\rho c_1)t}{\sqrt{1-\rho^2} }} \sim \Phi\left( \frac{\rho c_1-c_2 }{\sqrt{1-\rho^2}}\right).$$
With \cite{DHK20} [Thm 2.1] we have that
$$\pi_{[0,1]^2 \setminus \overline{\Delta}(u)^2,\rho}(c_1,c_2;u, au)=o\left(\pi_{[1-\frac{\Delta}{u^2},1],\rho}(c_1,u)\right)$$
and therefore asymptotics follows from \eqref{1dressoj}.
\QED
\subsection{Proof of Theorem \ref{MainSojurn}}
Let
$$\Sigma_{ s,t}= \begin{pmatrix}
s&  \rho \min(s,t) \\
\rho \min(s,t) &  t
\end{pmatrix} $$
be the covariance matrix of $(W_1(s),W_2(t)).$  We introduce the optimization problem that was used in \cite{DHK20,krystecki2021parisian}. For $\vk a = (1+\frac{c_1 s}{u},a+\frac{c_2t}{u})^\top$ denote
 $$ q_{\vk a}(s,t):= \vk a ^\top \Sigma^{-1}_{s,t} \vk a ,\quad \vk{b}(s,t):=  \Sigma^{-1}_{s,t} \vk a$$
 and set
\bqn{
	 q_{\vk a}^*(s,t)= \min_{ \vk x \ge \vk a} q_{\vk x}(s,t), \quad  q_{\vk a}^*= \min_{s,t \in [0,1]}  q_{\vk a}^*(s,t).
	}
We have for $a > \rho$ and $u$ large enough that $\vk b(s,t) \sim (\frac{t-a\rho\min(s,t)}{st-\rho^2(\min(s,t))^2},\frac{as-\rho\min(s,t)}{st-\rho^2(\min(s,t))^2}) > \vk 0.$
\cite{DEBICKIKOSINSKI} gives us that for any $s,t$ positive
\bqn{ \label{EXA2}
 \limit{u} \frac{1}{u^2} 	\log \pk{\exists_{s,t \in [0,1]} W_1^*(s)> u, W_2^*(t) > au} =  - \frac{q_{\vk a}^*(s,t)}{2} .
}	
Behaviour of this function will reflect the asymptotics for the cumulative Parisian ruin. We further denote $$t^*:=\lim_{u \to \IF}t_u,$$
where $t_u$ is the solution of the optimization problem of $q_{\vk a}^*(s,t)$ solved in \cite{DHK20}[Lemma 3.1].
As in the proof of Theorem \ref{simpleSojurn} we again can focus on the asymptotics of $\mathcal{S}_{S_1,S_2}(c_1,c_2,u,au).$ The following lemma tackles the behaviour of the probability of ruin on the square areas of size $\Theta(\frac{\Delta}{u^2})$ for $\Delta>0.$ Let us introduce a shorter notation $\eta_{u,k,l}(s,t):=(\eta_{1,u,k}(s),\eta_{2,u,l}(t)):=u(W_1(\frac{s}{u^2}+k_u) - W_1(k_u) - c_1\frac{s}{u^2}, W_2(\frac{t}{u^2}+l_u )- W_2(l_u) - c_2\frac{t}{u^2}).$
\BEL \label{PickandsSojurn}
Let $\rho \in (-1,1), a \in (\max(0,\rho),1], l,k=O(\frac{u\log(u)}{\Delta}) $ and  $\Delta, S_1,S_2>0 $ be given constants. Set $\mu_u=u^{-2} \varphi_{t_u}(u+c_1 ,au+c_2 t_u)$. Then, as $u\to \IF$
\BQNY
\mathcal{S}_{E_{u,k,l},\frac{(S_1,S_2)}{u^2}}(c_1,c_2,u,au)&\sim &\mu_u I_2(\Delta)  e^{-\frac{1}{2}u^2 (q_a(k_u,l_u)-q_a(1,t_u))},
\EQNY
where

$I_2(\Delta)= 	\begin{cases}
\int_{\R^2}\pk*{\int_{[0,\Delta]}\mathbf{1}(W_1(s)-s>x)ds>S_1 ,\int_{[0,\Delta]}\mathbf{1}(W_2(t)-at> y)dt>S_2}  e^{\lambda_1 x + \lambda_2 y} dxdy &l_u = k_u\\
\begin{aligned}\int_{\R^2}\pk*{\int_{[0,\Delta]}\mathbf{1}(W_1(s)-s>x)ds>S_1 }\\
\times \pk*{\int_{[0,\Delta]}\mathbf{1}( W_2(t)-\frac{a-\rho}{t^*-\rho^2}t> y)dt>S_2 }  e^{\lambda_1 x + \lambda_2 y} dxdy \end{aligned} &l_u > k_u \\
\begin{aligned}\int_{\R^2}\pk*{\int_{[0,\Delta]}\mathbf{1}(W_1(s)-\frac{1-a\rho}{1-\rho^2t^*}s>x)ds>S_1 }\\
\times \pk*{\int_{[0,\Delta]}\mathbf{1}( W_2(t)-\frac{a}{t^*}t> y)dt>S_2 }  e^{\lambda_1 x + \lambda_2 y} dxdy \end{aligned} &l_u < k_u
\end{cases}$ \\
and $\lambda_1 = \begin{cases}
\frac{1}{t^*}\frac{1-a\rho}{1-\rho^2} &l_u = k_u\\
\frac{t^*-a\rho}{t^*-\rho^2}, &l_u > k_u \\
\frac{1-a\rho}{1-\rho^2t^*}, &l_u < k_u \\
\end{cases}$, $\lambda_2 = \begin{cases}
\frac{1}{t^*}\frac{a-\rho}{1-\rho^2}, &l_u = k_u\\
\frac{a-\rho}{t^*-\rho^2 }, &l_u > k_u \\
\frac{a-\rho t^*}{t^* -\rho^2(t^*)^2}, &l_u < k_u \\
\end{cases}.$

Moreover
\begin{equation}
\label{Int2}
 \lim_{u \to \infty}\sup_{l,k = O(u\log u)}\int_{\R^2}\pk*{\begin{array}{ccc}
	\int_{E}\mathbf{1}(\eta_{1,u,k}(s)>x)ds>S_1  \\
	\int_{E}\mathbf{1}(\eta_{2,u,l}(t)> y)dt>S_2
	\end{array}
    \Bigg{|}
    \begin{array}{ccc}
	W_1^*(k_u )= u - \frac{x}{u}  \\
	W_2^*(l_u )= au - \frac{y}{u}
	\end{array} }  e^{\lambda_1 x + \lambda_2 y} dxdy\\
< \infty.
\end{equation}
\EEL
\prooflem{PickandsSojurn}
We can write
\BQNY
\lefteqn{\mathcal{S}_{E_{u,k,l},\frac{(S_1,S_2)}{u^2}}(c_1,c_2,u,au)}\\
&=&
\int_{\R^2}\pk*{\int_{E}\mathbf{1}(W_1^*(\frac{s}{u^2}+ k_u )>u)ds>S_1 ,\int_{E}\mathbf{1}(W_2^*(\frac{t}{u^2}+ l_u )> au)dt>S_2
    \Bigg{|}
    \begin{array}{ccc}
	W_1^*(k_u )= u - \frac{x}{u}  \\
	W_2^*(l_u )= au - \frac{y}{u}
	\end{array} } \\	
&\times& u^{-2}\varphi_{k_u,l_u}( u + c_1k_u - \frac{x}{u} ,  au + c_2l_u - \frac{y}{u} ) dxdy\\
&=&
\int_{\R^2}\mathbb{P}\Bigg{\{}\int_{E}\mathbf{1}(W_1^*(\frac{s}{u^2}+k_u ) - W_1(k_u ) +W_1(k_u)>u)ds>S_1 ,\\
& &\int_{E}\mathbf{1}(W_2^*(\frac{t}{u^2}+l_u )- W_2(l_u ) +W_2(l_u)> au)dt>S_2
    \Bigg{|}
    \begin{array}{ccc}
	W_1^*(k_u )= u - \frac{x}{u}  \\
	W_2^*(l_u )= au - \frac{y}{u}
	\end{array} \Bigg{\}} \\
&\times& u^{-2}\varphi_{k_u,l_u}( u + c_1k_u - \frac{x}{u} ,  au + c_2l_u - \frac{y}{u} ) dxdy\\
&=& \int_{\R^2}\mathbb{P}\Bigg{\{}\int_{E}\mathbf{1}(\eta_{1,u,k}(s)>u)ds>S_1 ,\int_{E}\mathbf{1}(\eta_{2,u,l}(t)> au)dt>S_2
    \Bigg{|}
    \begin{array}{ccc}
	W_1^*(k_u )= u - \frac{x}{u}  \\
	W_2^*(l_u )= au - \frac{y}{u}
	\end{array} \Bigg{\}} \\
&\times& u^{-2}\varphi_{k_u,l_u}( u + c_1k_u - \frac{x}{u} ,  au + c_2l_u - \frac{y}{u} ) dxdy
\EQNY
Notice that \cite{DHK20} [Lemma 3.3] gives us both for $k_u>l_u$ and $k_u<l_u,$ that as $u \to \IF$
\BQN \label{density}
\varphi_{k_u,l_u}( u + c_1 k_u - \frac{x}{u} ,  au + c_2 l_u - \frac{y}{u}    )&\sim& \varphi_{t_u}( u + c_1 ,  au + c_2 t_u) \\
&&\times e^{-\frac{1}{2} u^2(q_{\vk a_u(k_u,l_u)}(k_u,l_u)-q_{\vk a_u(1,t_u)}(1,t_u))}e^{\lambda_1 x + \lambda_2 y}. \nonumber
\EQN
Next we investigate the behaviour of the process
$$\Bigg{\{}\begin{array}{ccc}
	\eta_{1,u,k}(s)  \\
	\eta_{2,u,l}(t)
	\end{array} \Bigg{|}
    \begin{array}{ccc}
	W_1^*(k_u )= u - \frac{x}{u}  \\
	W_2^*(l_u )= au - \frac{y}{u}
	\end{array} \Bigg{\}}.$$
This process has been already studied in \cite{DHK20} [Lemma 3.3] wherein it distribution was determined to be Gaussian with parameters depending on whether $k_u<l_u, k_u=l_u$ or $k_u>l_u.$ However, to use the aforementioned calculations we need to prove that \eqref{Int2} is finite, so that we can use dominated convergence theorem. This comes straightforwardly from the combination of \eqref{comparison} and the fact that \cite{DHK20}(3.8) is finite. This completes the proof.
\QED

Before we move to the main proof, we introduce two lemmas that tackle the finitness of the constants that come from \nelem{PickandsSojurn}.
\BEL \label{1dFiniteSojurn}
i) For any $b, c>0, S \ge 0$ such that  $2b>c$  we have
$$
\int_{\R} \pk{ \int_{[0,T]} \mathbf{1}(W(t) - bt> x)dt>S } e^{ cx} dx \in (0,\infty).
$$
ii) For any $b>0, S \ge 0$
$$
\lim_{T\to\infty}\frac{1}{T}\int_{\R} \pk{ \int_{[0,T]} \mathbf{1}(W(t) - bt> x )dt>S} e^{ 2bx} dx \in (0,\infty).
$$
\EEL
Finitness and positivity of the one-dimensional constants has been proven in \cite{debicki2019sojourn}[proof of Thm 3.4].
\BEL \label {2dFiniteSojurn} Take any $a > \max(0,\rho), S_1,S_2 \ge 0.$ Then
$$
\int_{\R^2}\pk*{\int_{\R_+}\mathbf{1}(W_1(s)-s>x)ds>S_1 ,\int_{\R_+}\mathbf{1}(W_2(t)-at> y)dt>S_2}  e^{\lambda_1 x + \lambda_2 y} dxdy \in (0,\infty)
$$
where $\lambda_1 = \frac{1-a\rho}{1-\rho^2}, \lambda_2 = \frac{a-\rho}{1-\rho^2}.$
\EEL
Notice that the constants in \cite{DHK20} are analogous with the exchange of
$$\int_{A} \mathbf{1}(\vk{W}(\vk{t}) - \vk{b}\vk{t}> \vk{x})d\vk{t}>\vk{S}$$
with the classical ruin part
$$\sup_{\vk{t} \in A}\vk{W}(\vk{t}) - \vk{b}\vk{t}> \vk{x}$$
and hence finiteness of all the constants above follows from \eqref{comparison}. Positivity follows from the fact that cumulative Parisian ruin is not smaller than Parisian ruin and by positivity of the constants in \cite{krystecki2021parisian}, which are the same with the replacement of the cumulative Parisian ruin part with the Parisian ruin.

\prooftheo{MainSojurn}
Introduce
$$N_u:=\floor{\frac{u\log(u)}{\Delta}}, \quad E_{u,m}^1:=[(m+1)_u, m_u],\quad E_{u,j}^2:=[(j+1)_u, j_u],$$
$$K_u^{(1)}=\frac{(c_2-c_1\rho)u}{\Delta}, \quad K_u^{(2)}=\frac{(c_1-c_2\rho)u}{\Delta},$$
where $m_{u}=  1- \frac{(m-1)\Delta}{u^2}, j_{u}=  t^*- \frac{(j-1)\Delta}{u^2}.$ Different cases in the theorem are the result of various types of behaviour of function $q_{\vk a}^*(s,t)$ around the optimizing point.
\newline
\underline{Case (i): Suppose that $\rho>\frac{1}{4a}(1-\sqrt{8a^2+1}).$} According to \cite{DHK20}[Lemma 3.1] $t^*=1$. Denote $F_u:=E_{u,1}^2.$ For $\Delta>0$ we have
$$\mathcal{S}_{[0,1]^2,\frac{(S_1,S_2)}{u^2}}(c_1,c_2;u, au) \ge \mathcal{S}_{F_u,\frac{S}{u^2}}(c_1,c_2;u, au).$$
On the other hand
$$\mathcal{S}_{[0,1]^2,\frac{(S_1,S_2)}{u^2}}(c_1,c_2;u, au) \le \mathcal{S}_{F_u,\frac{S}{u^2}}(c_1,c_2;u, au)+\pi_{[0,1]^2 \setminus F_u}(c_1,c_2;u, au).$$
Using \nelem{PickandsSojurn} and \nelem{2dFiniteSojurn} and taking $u \to \infty$ and $\Delta \to \infty$, we get that
\BQNY
\mathcal{S}_{E_{u,1}^2,\frac{S}{u^2}}(c_1,c_2;u, au) &\sim &  \int_{\R^2}\pk*{\int_{[0,\infty)}\mathbf{1}(W_1(s)-s>x)ds>S ,\int_{[0,\IF)}\mathbf{1}(W_2(t)-at> y)dt>S}  \\
&\times& e^{\lambda_1 x + \lambda_2 y} dxdy u^{-2} \varphi_{1}(u+c_1 ,au+c_2).
\EQNY
From \cite{DHK20}, [Thm. 2.2] we have that
$$\pi_{[0,1]^2 \setminus E_{u,1}^2}(c_1,c_2;u, au)=o(u^{-2} \varphi_{1}(u+c_1 ,au+c_2) )=o(\mathcal{S}_{E_{u,1}^2,\frac{S}{u^2}}(c_1,c_2;u, au)).$$
With that, the proof of case (i) is complete.
\newline
\underline{Case (ii): Suppose that $\rho=\frac{1}{4a}(1-\sqrt{8a^2+1}).$} From this case onwards we observe a quadratic behaviour of the variance near the optimal point, hence more than one square area of size $O\left(\frac{\Delta}{u^2}\right)$ impacts the asymptotics. From case (ii) of \cite{krystecki2021parisian}[Thm. 2.2] recall that
$$u^2(q_a(k_u,l_u)-q_a(1,1))=\tau_1 (k-1)\Delta+\tau_4 \frac{(l-1)^2\Delta^2}{u^2}+o\left(\frac{k^2}{u^2}\right)+o\left(\frac{l^3}{u^4}\right),$$
where $\tau_1=\frac{(1 - a \rho)^2}{(1 - \rho^2)^2}>0$ and $\tau_4=\frac{\rho^2  - 2 a \rho^3 +  a^2 \rho^2 }{(1 - \rho^2)^2}>0.$ Constants $\tau_i$ will vary for the other cases, analogously to what was calculated in \cite{krystecki2021parisian}[Thm. 2.2]. This case is split into two subcases. First consider $c_2-\rho c_1 \le 0.$ According to \cite{DHK20}[Lemma 3.1] $t^*=1$. For $F_u:=[1-\frac{\Delta}{u^2},1]\times[1-\frac{\log(u)}{u},1-\frac{\Delta}{u^2}]$ we write
$$\mathcal{S}_{[0,1]^2,\frac{(S_1,S_2)}{u^2}}(c_1,c_2;u, au) \ge \mathcal{S}_{F_u,\frac{S}{u^2}}(c_1,c_2;u, au).$$
Since either the whole period of crossing the barrier occurs on $F_u$ or we cross the barrier at least once on $[0,1]^2 \setminus F_u,$ then
$$\mathcal{S}_{[0,1]^2,\frac{(S_1,S_2)}{u^2}}(c_1,c_2;u, au) \le \mathcal{S}_{F_u,\frac{S}{u^2}}(c_1,c_2;u, au)+\pi_{[0,1]^2 \setminus F_u}(c_1,c_2;u, au).$$
Using Bonferroni inequality we have that
\BQN \label{lBS}
\mathcal{S}_{F_u,\frac{(S_1,S_2)}{u^2}}(c_1,c_2;u, au) &\ge& \sum_{l=2}^{N_u}\pk*{\int_{E_{u,1}^1}\mathbf{1}(W_1^*(s)>u)ds>\frac{S_1}{u^2},\int_{E_{u,l}^2}\mathbf{1}(W_2^*(t)> au)dt>\frac{S_2}{u^2}} \nonumber \\
&&- \sum_{l=2}^{N_u}\sum_{m=l+1}^{N_u}\pk*{\exists_{s \in E_{u,1}^1, t_1 \in E_{u,l}^2, t_2 \in E_{u,m}^2}:\begin{array}{ccc} W_1^*(s)>u \\ W_2^*(t_1)> au \\ W_2^*(t_2)> au \end{array} } \nonumber\\
&:=&S_{u,\Delta}-D_{u,\Delta}
\EQN
Further we have
\BQN \label{uBS}
\mathcal{S}_{F_u,\frac{(S_1,S_2)}{u^2}}(c_1,c_2;u, au) &\le& S_{u,\Delta}+D_{u,\Delta}.
\EQN
From \nelem{PickandsSojurn} we have as $u \to \IF$
\BQNY
S_{u,\Delta}& \sim & C_{2,\mathcal{S}}^{(1)}(\Delta) C_{2,\mathcal{S}}^{(2)}(\Delta) u^{-2}\varphi_{t^*}(u+c_1 ,au+c_2) \sum_{l=2}^{N_u} e^{-\frac{1}{2}u^2(q_{\vk a}(k_u,l_u)-q_{\vk a}(1,1))},\\
\EQNY
where
$$C_{2,\mathcal{S}}^{(1)}(\Delta)=\int_{\R}\pk*{ \int_{[0,\Delta]}\mathbf{1}( W_1(s) - \frac{1-a\rho}{1-\rho^2}s> x)ds>S_1} e^{ \frac{1-a\rho}{1-\rho^2}x} dx$$
and
$$C_{2,\mathcal{S}}^{(2)}(\Delta)=\int_{\R} \pk*{\int_{[0,\Delta]}\mathbf{1}( W_2(t) - at> x)dt>S_2} e^{ 2ax} dx.$$
Using Taylor expansions together with \cite{DHK20}[Lem 3.6] we get that as $u \to \IF$
\BQNY
S_{u,\Delta} & \sim & C_{2,\mathcal{S}}^{(1)}(\Delta) C_{2,\mathcal{S}}^{(2)}(\Delta) u^{-2}\varphi_{t^*}(u+c_1 ,au+c_2) \sum_{l=2}^{N_u} e^{-\frac{\tau_4}{2}\frac{(l-1)^2\Delta^2}{u^2}}\\
&=& \frac{1}{\sqrt{\tau_4}}C_{2,\mathcal{S}}^{(1)}(\Delta) \frac{C_{2,\mathcal{S}}^{(2)}(\Delta)}{\Delta} u^{-1}\varphi_{t^*}(u+c_1 ,au+c_2) \sum_{l=2}^{N_u} \frac{\sqrt{\tau_4}\Delta}{u} e^{-\frac{\tau_4}{2}\frac{(l-1)^2\Delta^2}{u^2}}\\
&\sim& C_{2,\mathcal{S}}^{(1)}(\Delta) \frac{C_{2,\mathcal{S}}^{(2)}(\Delta)}{\Delta} \frac{\sqrt{\pi}}{\sqrt{2\tau_4}} u^{-1}\varphi_{t^*}(u+c_1 ,au+c_2).
\EQNY
With \nelem{1dFiniteSojurn} we have
$$\lim_{\Delta \to \IF}C_{2,\mathcal{S}}^{(1)}(\Delta)=C_{2,\mathcal{S}}^{(1)}(\IF) \in (0,\IF),\quad \lim_{\Delta \to \IF}\frac{C_{2,\mathcal{S}}^{(2)}(\Delta)}{\Delta}=C_{2,\mathcal{S}}^{(2)}(\IF) \in (0,\IF).$$
Hence
$$\lim_{\Delta \to \IF}\lim_{u \to \IF}\frac{S_{u,\Delta}}{C_{2,\mathcal{S}}^{(1)}(\IF)C_{2,\mathcal{S}}^{(2)}(\IF)\frac{\sqrt{\pi}}{\sqrt{2\tau_4}} u^{-1}\varphi_{t^*}(u+c_1 ,au+c_2)}=1.$$
From \cite{DHK20}[Theorem 2.2, case (ii)] we have that
\BQN \label{negl1}
\lim_{\Delta \to \IF}\lim_{u \to \IF}\frac{D_{u,\Delta}}{S_{u,\Delta}}&=&\lim_{\Delta \to \IF}\lim_{u \to \IF}\frac{D_{u,\Delta}}{C_{2,\mathcal{S}}^{(1)}(\IF) C_{2,\mathcal{S}}^{(2)}(\IF) \frac{\sqrt{\pi}}{\sqrt{2\tau_4}}u^{-1}\varphi_{t^*}(u+c_1 ,au+c_2)}=0
\EQN
and also that
\BQN \label{negl2}
\lim_{\Delta \to \IF}\lim_{u \to \IF}\frac{\pi_{[0,1]^2 \setminus F_u}(c_1,c_2;u, au)}{S_{u,\Delta}}&=&\lim_{\Delta \to \IF}\lim_{u \to \IF}\frac{\pi_{[0,1]^2 \setminus F_u}(c_1,c_2;u, au)}{C_{2,\mathcal{S}}^{(1)}(\IF) C_{2,\mathcal{S}}^{(2)}(\IF) \frac{\sqrt{\pi}}{\sqrt{2\tau_4}}u^{-1}\varphi_{t^*}(u+c_1 ,au+c_2)}=0.
\EQN
Therefore, using the asymptotics for $\pi_{[0,1]^2}(c_1,c_2;u, au)$ from \cite{DHK20}[Thm 2.2] we have
$$\lim_{u \to \IF}\mathscr{S}_{S_1,S_2}(c_1,c_2;u,au) = \frac{(1-a\rho)\widehat{\mathcal{P}}(\frac{1-a\rho}{1-\rho^2},\frac{1-a\rho}{1-\rho^2},S_1)\widehat{\mathcal{H}}(a,2a,S_2)}{2a(1-\rho^2)}.$$
It remains to consider the case $c_2-\rho c_1 > 0.$ Using \cite{DHK20}[Lemma 3.1] we have that the only minimizer of $q_{\vk a_u(s,t)}^*(s,t)$  on $[0,1]^2$ is $(s_u,t_u)=\left(1,\frac{a}{\rho(2a\rho-1)+\frac{c_2-\rho c_1}{u}}\right)$ with
$$\frac{a}{\rho(2a\rho-1)+\frac{c_2-\rho c_1}{u}} \nearrow 1$$
as $u \to \IF.$ The optimal area according to \cite{DHK20}[Thm 2.2, case (iii)] is $F_u:=[1-\frac{\Delta}{u^2},1]\times[t_u-\frac{\log(u)}{u},1-\frac{\Delta}{u^2}].$
From \nelem{PickandsSojurn} we have as $u \to \IF$
\BQNY
S_{u,\Delta}& \sim & C_{2,\mathcal{S}}^{(1)}(\Delta) C_{2,\mathcal{S}}^{(2)}(\Delta) u^{-2}\varphi_{t^*}(u+c_1 ,au+c_2) \sum_{l=-K_u^{(1)}}^{N_u} e^{-\frac{1}{2}u^2(q_{\vk a}(k_u,l_u)-q_{\vk a}(1,1))},\\
\EQNY
where $C_{2,\mathcal{S}}^{(1)}(\Delta),C_{2,\mathcal{S}}^{(2)}(\Delta)$ are defined as above. Using Taylor expansions and \cite{DHK20}[Lem 3.6] we get that as $u \to \IF$
\BQNY
S_{u,\Delta} & \sim & C_{3,\mathcal{S}}^{(1)}(\Delta) C_{3,\mathcal{S}}^{(2)}(\Delta) u^{-2}\varphi_{t^*}(u+c_1 ,au+c_2) \sum_{l=-K_u^{(1)}}^{N_u} e^{-\frac{\tau_4}{2}\frac{(l-1)^2\Delta^2}{u^2}}\\
&\sim& C_{2,\mathcal{S}}^{(1)}(\Delta) \frac{C_{2,\mathcal{S}}^{(2)}(\Delta)}{\Delta} \frac{\sqrt{2\pi}}{\sqrt{\tau_4}} e^{-a\frac{(c_1\rho-c_2)^2}{2\rho(1-a\rho)}}  \Phi\left(c_2-\rho c_1\right) u^{-1}\varphi_{t^*}(u+c_1 ,au+c_2).
\EQNY
Using \eqref{uBS} and \eqref{lBS} together with \cite{DHK20}[Theorem 2.2, case (iii)] we have that
$$\lim_{\Delta \to \IF}\lim_{u \to \IF}\frac{S_{u,\Delta}}{C_{2,\mathcal{S}}^{(1)}C_{2,\mathcal{S}}^{(2)}\frac{\sqrt{2\pi}}{\sqrt{\tau_4}} e^{-a\frac{(c_1\rho-c_2)^2}{2\rho(1-a\rho)}}  \Phi\left(c_2-\rho c_1\right)  u^{-1}\varphi_{t^*}(u+c_1 ,au+c_2)}=1$$
and that \eqref{negl1} and \eqref{negl2} hold. Together with the asymptotics for $\pi_{[0,1]^2}(c_1,c_2;u, au)$ from \cite{DHK20}[Thm 2.2] we have
$$\lim_{u \to \IF}\mathscr{S}_{S_1,S_2}(c_1,c_2;u,au) = \frac{(1-a\rho)\widehat{\mathcal{P}}(\frac{1-a\rho}{1-\rho^2},\frac{1-a\rho}{1-\rho^2},S_1)\widehat{\mathcal{H}}(a,2a,S_2)}{2a(1-\rho^2)},$$
which completes the proof of case (ii). \newline
The following cases follow the same path of proof, where the varying component is the main area $F_u$ that we use in analogons of \eqref{uBS} and \eqref{lBS}. Hence the proofs for the remaining cases are omitted apart from pointing out the main area and the optimizing point for the function $q_{\vk a_u(s,t)}^*(s,t)$ on $[0,1]^2.$ \newline
\underline{Case (iii): Suppose that $\rho=-\frac{1}{2}, a=1.$} According to \cite{DHK20}[Lem. 3.1] $t^*=1$. From the proof of \cite{DHK20}[Thm 2.2 case (iv)] we have that
$$F_u:=[1-\frac{\Delta}{u^2},1]\times[1-\frac{\log(u)}{u},1-\frac{\Delta}{u^2}]\cup [1-\frac{\log(u)}{u},1-\frac{\Delta}{u^2}]\times[1-\frac{\Delta}{u^2},1].$$
\underline{Case (iv): Suppose that $\rho<\frac{1}{4a}(1-\sqrt{8a^2+1}).$} \cite{DHK20}[Lem. 3.1] gives exactly one minimizer of
$q_{\vk a_u(s,t)}^*(s,t)$  on $[0,1]^2$ which is $(s_u,t_u)=(1,\frac{a}{\rho(2a\rho-1)+\frac{c_2-\rho c_1}{u}})$ and as $u \to \IF$ we have $t_u<1.$
Proof of \cite{DHK20}[Thm 2.2 case (v)] leads to the main area being
$$F_u:=[1-\frac{\Delta}{u^2},1]\times [t_u-\frac{\log(u)}{u},t_u+\frac{\log(u)}{u}].$$
\underline{Case (v): Suppose that $a=1, \rho< A_a.$} From \cite{DHK20}[Lem. 3.1] there are two optimal points:
$$(s_u,t_u)= (1,\frac{1}{\rho(2\rho-1)+\frac{c_2-\rho c_1}{u}}), \quad (\bar{s}_u,\bar{t}_u)= (\frac{1}{\rho(2\rho-1)+\frac{c_1-\rho c_2}{u}},1).$$
The main area is a combination of two disjoint areas.
$$F_u:=[1-\frac{\Delta}{u^2},1]\times[t_u-\frac{\log(u)}{u},t_u+\frac{\log(u)}{u}]\cup [t_u-\frac{\log(u)}{u},t_u+\frac{\log(u)}{u}]\times[1-\frac{\Delta}{u^2},1].$$
This completes the proof.
\QED

\bibliographystyle{ieeetr}

\bibliography{queue2d}
\end{document}